\documentclass[a4paper, 12pt]{article}
\usepackage{enumerate, theorem}
\usepackage{amsmath, amsfonts, amssymb}
\usepackage[height=22.5cm, width=15cm]{geometry}
\usepackage[all]{xy}
\usepackage{mathrsfs, yfonts}

\newcommand{\parag}[1]{\paragraph{\sc{#1.}} }
\def\A{{\tilde{\mathcal{A}}}}

\newtheorem{thm}{Th\'eor\`eme}[subsection] 
\newtheorem{defn}[thm]{D\'efinition}

\newtheorem{prop}[thm]{Proposition}
\newtheorem{lemma}[thm]{Lemme}

\begin{document}

\title{Sur les fonctions a singularit\'e de dimension 1.}

\author{Daniel Barlet\footnote{Barlet Daniel, Institut Elie Cartan UMR 7502  \newline
Nancy-Universit\'e, CNRS, INRIA  et  Institut Universitaire de France, \newline
BP 239 - F - 54506 Vandoeuvre-l\`es-Nancy Cedex.France.\newline
e-mail : barlet@iecn.u-nancy.fr}.}

\date{version revis\'ee le  03/09/07 \\(seconde version de la pr\'epublication 2007/01)}

\maketitle

\markright{Sur les fonctions a singularit\'e de dimension 1.}

\section*{Abstract}
In this article we show that all results proved for a large class of holomorphic germs \ $f : (\mathbb{C}^{n+1}, 0) \to (\mathbb{C}, 0)$ \ with a 1-dimension singularity in [B.II]  are valid for an arbitrary such germ. 

\bigskip

AMS Classification (2000) : 32-S-25, 32-S-40, 32-S-50.

\bigskip

\tableofcontents

\section{Introduction.}

\subsection{G\'en\'eralisation de l'\'etude locale aux points g\'en\'eriques de la courbe  \ $S$.}

Dans le pr\'esent article nous consid\`ererons la situation suivante : \\
Soit \ $ f : ( \mathbb{C}^{n+1}, 0) \to (\mathbb{C}, 0) $ \ un germe de fonction holomorphe dont le lieu singulier \ $S : = \{ x \ / \ f(x) = 0 \ {\rm et} \ df_x = 0 \}$ \ est un germe de courbe \`a l'origine, et  soit  \ $t : ( \mathbb{C}^{n+1}, 0) \to (\mathbb{C}, 0) $ \ un germe de fonction  holomorphe non singuli\`ere tel que la restriction de \ $t$ \ \`a \ $(S, 0)$ \ soit  finie. Sur un voisinage ouvert \ $X$ \ de l'origine assez petit on aura pour chaque \ $\sigma \in S^* : = S \setminus \{0\}$ \ un (a,b)-module associ\'e au germe en \ $\sigma$ \  de fonction  \`a singularit\'e isol\'ee obtenu en restreignant \ $f$ \ \`a l'hypersurface lisse d\'efinie par \ $\{ t = t(\sigma)\}$. Rappelons que ce (a,b)-module n'est rien d'autre que le compl\'et\'e formel du module de Brieskorn de ce germe de fonction \`a  singularit\'e isol\'ee (voir par exemple [B.05]).\\
Les faisceaux de cohomologie \ $\hat{\mathcal{H}}^{\bullet}$ \  du complexe de de Rham
$$ (\hat{K}er\, df^{\bullet}, d^{\bullet}) $$
o\`u \ $\hat{K}er\, df^{\bullet}$ \ d\'esigne le noyau de la multiplication ext\'erieure par \ $df$
$$ \wedge df^{\bullet} : \hat{\Omega}^{\bullet} \to  \hat{\Omega}^{\bullet +1}$$
agissant sur le compl\'et\'e formel en \ $f$ \ des formes holomorphes, sont nuls en degr\'es diff\'erents de  \ $1, n, n+1$ \ dans cette situation et les faisceaux \ $\hat{\mathcal{H}}^n$ \ et \ $\hat{\mathcal{H}}^{n+1}$ \ sont \`a support dans \ $S$. La consid\'eration du germe auxilliaire \ $t$ \ permet de calculer  ces deux faisceaux comme respectivement  noyau et conoyau d'une $t-$connexion compatible avec les op\'erations \  $a$ \ et  \ $b$ \ de la fa{\c c}on suivante.

\smallskip

Soit \ $(\hat{K}er\, {df_/}^{\bullet}, {d_/}^{\bullet})$ \ le complexe de de Rham $t-$relatif compl\'et\'e formellement en  $f$, et soit \ $\mathbb{E} : = \mathcal{H}^n(\hat{K}er\, {df_/}^{\bullet}_/, {d_/}^{\bullet})$. Ce faisceau est naturellement muni d'op\'erations $a$ \ et \ $b$ \ ainsi que d'une structure de \ $t^{-1}(\mathcal{O}_D)-$module, o\`u \ $D$ \ est un disque de centre \ $0$ \ assez petit dans \ $\mathbb{C}$. Nous d\'efinissons alors  une \ $t^{-1}(\mathcal{O}_D)-$connexion
$$ b^{-1}.\nabla :  \mathbb{P} \to  \mathbb{E} $$
o\`u \ $\mathbb{P} : = \{ x \in \mathbb{E}\ / \  \nabla(x) \in b.\mathbb{E}\}$, qui commute \`a \ $a$ \ et \ $b$. \\
La proposition suivante est d\'emontr\'ee dans [B.II] prop. 4.2.8.

\begin{prop}\label{connexion}
On a une suite exacte naturelle de faisceaux de \ $\underline{\mathbb{C}}_X-$modules \`a support \ $S$ :
$$ 0 \to \hat{\mathcal{H}}^n \overset{u}{\to} \mathbb{P} \overset{b^{-1}.\nabla}{\longrightarrow} \mathbb{E} \overset{v}{\to} \hat{\mathcal{H}}^{n+1} \to 0 $$
compatible aux op\'erations \ $a$ \ et \ $b$ \ sur ces faisceaux, o\`u \ $u$ \ est induite par la projection \'evidente \ $ \hat{K}er\, df^n \to \hat{\Omega}_/^n $ \ et \ $v$ \ par la multiplication ext\'erieure par \ $dt$.
\end{prop}

L'objectif principal de cet article est de montrer le r\'esultat "technique"  suivant :

\begin{thm}\label{technique}
Dans la situation pr\'ecis\'ee ci-dessus, on a
\begin{enumerate}
\item le faisceau \ $\hat{\mathcal{H}}^n$ \ est un syst\`eme local de (a,b)-modules g\'eom\'etriques sur \ $S^*$.
\item La restriction \`a \ $S^*$ \ du faisceau \ $\hat{\mathcal{H}}^{n+1}$ \ v\'erifie la propri\'et\'e de prolongement analytique suivante\\

(PA) $\quad$  Soient \ $V \subset U$ \ deux ouverts de \ $S^*$, avec \ $U$ \ connexe et \ $V$ \ non vide. Alors la restriction \ $ \Gamma(U, \hat{\mathcal{H}}^{n+1}) \to \Gamma(V, \hat{\mathcal{H}}^{n+1})$ \ est injective.
\end{enumerate}
\end{thm}

Ce r\'esultat montre que l'assertion " le faisceau \ $\hat{\mathcal{H}}^n$ \ est un syst\`eme local de (a,b)-modules r\'eguliers et g\'eom\'etriques sur \ $S^*$ \ " du  th\'eor\`eme 4.3.1 ainsi que le th\'eor\`eme  4.3.2  de [B.II]  sont valables sans restriction pour un germe de fonction holomorphe \`a singularit\'e de dimension 1. Ceci implique que l'ensemble des r\'esultats de [B.II]  est valable dans ce cadre. En particulier le th\'eor\`eme de finitude 5.2.1, le corollaire 5.2.3 ainsi que les th\'eor\`emes 6.2.1 et 6.4.1.\\
Les \'enonc\'es g\'en\'eralisant les th\'eor\`emes 5.2.1 et 6.2.1 sont d\'etaill\'es  ci-apr\`es pour la commodit\'e du lecteur.

\subsection{Quelques cons\'equences.}

\subsubsection{Rappels.}

{\bf Dans ce qui suit \ $f$ \ d\'esignera un germe \`a l'origine de \ $\mathbb{C}^{n+1}$ \ de fonction holomorphe qui est suppos\'e r\'eduit.}\\

\noindent Pour \ $f$ \ r\'eduite, on a un isomorphisme naturel compatible \`a \ $a$ \ et \ $b$, $$ E_1 \otimes \mathbb{C}_Y \simeq \hat{\mathcal{O}}_X.df \cap Ker\, d \simeq \mathcal{H}^1\big((\hat{K}er\,df)^{\bullet}, d^{\bullet})\big) $$
o\`u l'on a pos\'e \ $E_1 : = \mathbb{C}[[z]].dz$, muni des op\'erations \ $a : = \times z$ \ et \ $b : = (\int_0^z ).dz$.\\
Plac\'e en degr\'e 1  ce faisceau d\'efini un sous-complexe, que nous noterons \ $ E_1 \otimes \mathbb{C}_Y[1]$,  du complexe  \ $\big((\hat{K}er\,df)^{\bullet}, d^{\bullet}\big)$.

\noindent Nous d\'efinirons le complexe \ $\big((\tilde{K}er\,df)^{\bullet}, d^{\bullet}\big)$ comme le quotient
$$ \big((\tilde{K}er\,df)^{\bullet}, d^{\bullet}\big) : = \big((\hat{K}er\,df)^{\bullet}, d^{\bullet}\big)\big/ E_1 \otimes \mathbb{C}_Y[1].$$

\noindent Notons  \ $\widehat{\mathcal{A}}$ \ la \ $\mathbb{C}-$alg\`ebre
\ $ \widehat{\mathcal{A}} : = \lbrace \sum_{\nu =0}^{\infty} P_{\nu}(a).b^{\nu} , \ {\rm avec} \  P_{\nu} \in \mathbb{C}[x] \rbrace$ \
dont le produit est d\'efini par les deux conditions suivantes :
\begin{itemize}
\item{1)} On a la relation de commutation \ $ a.b - b.a = b^2$.
\item{2)} La multiplication par \ $a$ \ (\`a gauche ou \`a droite) est continue pour la filtration \ $b-$adique.
\end{itemize}

On a alors le r\'esultat g\'en\'eral suivant (voir [B.II] th. 2.1.1) :

  \begin{thm}\label{Structure}
  Les complexes \ $\big((\tilde{K}er\,df)^{\bullet}, d^{\bullet}\big)$ \ et \ $\big((\hat{K}er\,df)^{\bullet}, d^{\bullet}\big)$ \ sont  canoniquement quasi-isomorphes \`a des complexes de faisceaux de  \ $\widehat{\mathcal{A}}-$modules (\`a gauche) sur \ $Y$ \ ayant des diff\'erentielles \ $\widehat{\mathcal{A}}-$lin\'eaires. \\
  De plus, via ces quasi-isomorphismes, la suite exacte
  $$ 0 \to  E_1 \otimes \mathbb{C}_Y[1] \to \big((\hat{K}er\,df)^{\bullet}, d^{\bullet}\big) \to \big((\tilde{K}er\,df)^{\bullet}, d^{\bullet}\big) \to 0 $$
  correspond \`a une suite exacte de complexes de faisceaux de \ $\widehat{\mathcal{A}}-$modules.
  \end{thm}
  
  \bigskip
  
    \noindent Ce th\'eor\`eme montre l'existence d'une structure naturelle de  \ $\widehat{\mathcal{A}}-$modules sur tout groupe d'hypercohomologie de ces deux complexes. Il donne \'egalement la \ $\widehat{\mathcal{A}}-$lin\'earit\'e (\`a gauche) des applications naturelles entre ces groupes.
    
   \subsubsection{Les r\'esultats principaux.}
    
  {\bf   Supposons maintenant que le lieu singulier \ $S$ \ de \ $f$ \ est une courbe.}\\
  
  \smallskip
  
  Rappelons qu'un \ $\A-$module \ $E$ \  est dit r\'egulier g\'eom\'etrique s'il est de type fini sur la sous-alg\`ebre \ $\mathbb{C}[[b]]$ \ de \ $\A$, et si son quotient par sa \ $b-$torsion (qui toujours est un (a,b)-module, puisque \ $\mathbb{C}[[b]]-$libre de type fini) est un (a,b)-module r\'egulier g\'eom\'etrique.\\
  Rappelons encore qu'un (a,b)-module \ $E$ \ est r\'egulier si son satur\'e par \ $b^{-1}.a$, not\'e \ $E^{\sharp}$, est encore de type fini sur \ $\mathbb{C}[[b]]$. Il sera dit g\'eom\'etrique si, de plus, les valeurs propres de \ $b^{-1}.a$ \ agissant sur l'espace vectoriel de dimension finie \ $E^{\sharp}/b.E^{\sharp}$ \ sont dans \ $\mathbb{Q}^{+*}-1$ \ (donc rationnelles strictement sup\'erieures \`a \ $-1$). Cette condition est v\'erifi\'ee pour le (a,b)-module de Brieskorn d'un germe \`a singularit\'e isol\'ee gr\^ace au th\'eor\`eme de Monodromie et au th\'eor\`eme de positivit\'e de Malgrange.
    
    \bigskip

Notons  \ $ \hat{\mathcal{H}}^i$ \ le i-\`eme faisceau de cohomologie du complexe \ $\big((\hat{\mathcal{K}}er\,df)^{\bullet}, d^{\bullet}\big)$.\\
   Nous noterons \ $c \cap S$ \ (resp. \  $c$) \ la famille des ferm\'es de \ $Y : = f^{-1}(0)$ \ qui rencontrent \ $S$ \ suivant un compact (resp. la famille des compacts de \ $Y$).\\
    Le th\'eor\`eme 1.2.2 implique  la g\'en\'eralisation suivante du th\'eor\`eme 5.2.1  de [B.II].
    
   \begin{thm}\label{reg.}
   Les \ $\A-$modules suivants sont r\'eguliers g\'eom\'etriques :
   \begin{enumerate}
    \item \ $ H^i_{\{0\}}(Y, \hat{\mathcal{H}}^j)$ \ pour \ $ i \geq 0$ \ et \ $j = 1, n$.
    \item \ $E : = H^0_{\{0\}}(Y, \hat{\mathcal{H}}^{n+1})$.
     \item \quad $ E_{\Phi} : = \mathbb{H}_{\Phi}^{n+1}(Y, (\hat{K}er\,df)^{\bullet}, d^{\bullet}\big))$ \ et \quad \  $E'_{\Phi} : = \mathbb{H}_{\Phi}^{n+1}(Y, (\tilde{K}er\,df)^{\bullet}, d^{\bullet}\big))$ \\
    pour \ $\Phi = c$ \ et \ $\Phi = c\, \cap \, S $.
    \end{enumerate}
    \end{thm}
    
    On en d\'eduit  \'egalement la g\'en\'eralisation suivante du th\'eor\`eme 6.2.1 de [B.II].
    
    \begin{thm}\label{def. h}
On a pour \ $n \geq 2$ \ un accouplement (a,b)-sesquilin\'eaire naturel, non d\'eg\'en\'er\'e 
$$h : E'_{c\, \cap S} \times E \longrightarrow  \vert \Xi' \vert^2 $$
donn\'e par int\'egration dans les fibres de \ $f$ .
\end{thm}

On a pos\'e ici, pour \ $n \in \mathbb{N}$ \ fix\'e
\begin{equation*}
\vert \Xi' \vert^2 : = \sum_{r \in ]-1,0] \cap \mathbb{Q}\ , \ j \in [0,n]} \mathbb{C}[[s,\bar{s}]].\vert s \vert^{2r}.(Log \vert s \vert^2)^j \Big/ \mathbb{C}[[s,\bar{s}]]  \ .
\end{equation*}

\bigskip

\noindent Pr\'ecisons ce que signifie "donn\'e par int\'egration sur les fibres de \ $f$ ". \\
Soient \ $\omega \ {\rm resp.} \ \omega' $ \ des $(n+1)-$formes \ $ \mathcal{C}^{\infty} $ \ annul\'ees par \ $ \wedge df $ \ et par \ $d$, le support de \ $\omega'$ \ rencontrant \ $S$ \ suivant un compact \ $K$. Soit alors \ $ \rho \in  \mathcal{C}^{\infty}_c(X) $ \ v\'erifiant \ $ \rho \equiv 1 $ \ au voisinage de \ $K$ . Alors nous d\'efinirons
\begin{equation*}
 h([\omega'] , [\omega]) = \frac{1}{(2i\pi)^n} \int_{f=s}\rho. \frac{\omega'}{df}\wedge \overline{\frac{\omega}{df}} \  \in \vert \Xi' \vert^2 .
 \end{equation*}
 
 Le calcul des \ $\hat{\mathcal{A}}-$modules \ $E'_{c \cap S}$ \ et \ $E$ \ via de telles formes \ $\mathcal{C}^{\infty}$ \ est justifi\'e par le lemme 6.1.1 de [B.II].

\smallskip

\noindent Pr\'ecisons  enfin ce que nous entendons par accouplement ''non d\'eg\'en\'er\'e''.
\begin{itemize}
\item Pour tout \'element $ [\omega'] \in E'_{c\cap S} $ \ qui n'est pas de $b-$torsion, il existe un \'el\'ement  \ $[\omega]\in E \ {\rm telle \ que \ } h([\omega'],[\omega]) \not= 0  \ {\rm dans} \ \vert \Xi' \vert^2  .$
\item  Pour tout \'el\'ement  $ [\omega] \in E $ \ qui n'est pas de $b-$torsion, il existe un \'el\'ement \ $[\omega] \in E'_{c\cap S} \ {\rm tel \ que \ } h([\omega'],[\omega]) \not= 0  \ {\rm dans} \ \vert \Xi' \vert^2  .$
\end{itemize}

\section{La $b^{-1}.\nabla$ finitude.}

\subsection{La situation consid\'er\'ee} 
Soit  \ $f :  ( \mathbb{C}^{n+1}, 0) \to (\mathbb{C}, 0) $ \ un germe de fonction holomorphe \`a lieu singulier \ $(S,0)$ \ de dimension 1, et soit \ $ t : ( \mathbb{C}^{n+1}, 0) \to (\mathbb{C}, 0) $ \ un germe de fonction  holomorphe non singuli\`ere. Notons \ $(\Sigma, 0)$ \ le lieu critique du germe d'application $$(f,t) : ( \mathbb{C}^{n+1}, 0) \to (\mathbb{C}^2, 0).$$ 
Faisons l'hypoth\`ese que \ $(\Sigma, 0)$ \ soit de dimension 1 et que  la restriction de \ $t$ \ \`a \ $(\Sigma, 0)$ \ soit  finie.\\
Alors il existe un voisinage ouvert de Stein  \ $X$ \ de l'origine dans \ $\mathbb{C}^{n+1}$ \ et un disque \ $D$ \ de centre \ $0$ \ dans \ $\mathbb{C}$ \ v\'erifiant les propri\'et\'es suivantes :
\begin{enumerate}
\item Chaque composante irr\'eductible de \ $S : = \{ x \in X / \ f(x) = 0 \ {\rm et} \ df_x = 0 \}$ \ contient l'origine, est non singuli\`ere en dehors de l'origine et est un disque topologique.
\item La restriction \ $ t_{\vert \Sigma} : \Sigma  \to D$ \ est un rev\^etement ramifi\'e fini dont la restriction \ $  t_{\vert \Sigma^*} : \Sigma^* \to D^*$ \ est un rev\^etement non ramifi\'e, o\`u nous avons  pos\'e \\ $\Sigma^* : = \Sigma \setminus \{0\}$ \ et \ $D^* : = D \setminus \{0\}$.
\end{enumerate}

L'existence de tels voisinages ouverts \ $X$ \ et \ $D$ \ de \ $0$ \ dans \ $\mathbb{C}^{n+1}$ \ et \ $\mathbb{C}$ \ arbitrairement petits est \'el\'ementaire. L'existence, pour un germe \ $f$ \ donn\'e, de germes holomorphes auxiliaires convenables est cons\'equence du lemme suivant

 \begin{lemma}
 Soit \ $f$ \ un germe de fonction holomorphe \`a l'origine de \ $\mathbb{C}^{n+1}$ \ dont le lieu singulier \ $S$ \ est de dimension 1. Pour \ $l \in (\mathbb{C}^{n+1})^*$ \ dans un ouvert dense, le lieu critique \ $\Sigma_l$ \ du germe d'application \ $(f,l) :(\mathbb{C}^{n+1},0) \to (\mathbb{C}^2,0)$ \ est de dimension 1 et contient \ $S$;  on a donc l'\'egalit\'e \ $\Sigma_l = S$ \ au voisinage de \ $S^*$.\\
  De plus, on peut demander que  la restriction de \ $l$ \ \`a \ $\Sigma_l$ \ soit propre et finie avec un unique point de ramification \`a l'origine. \\
Un tel \ $l$ \ \'etant fix\'e, il existe  un disque ouvert \ $D$ \ centr\'e \`a l'origine dans \ $\mathbb{C}$ \ assez petit, tel que la famille \ $(f_{\tau})_{\tau \in D}$, o\`u \ $f_{\tau} : = f_{\vert \{ l = \tau \}}$,  soit sur \ $D^*$ \ une famille \`a \ $\mu-$constant le long de chaque composante connexe de \ $S^*$.
\end{lemma}

\parag{Preuve} Comme \ $(S,0)$ \ est par hypoth\`ese un germe de courbe, pour  \ $l \in (\mathbb{C}^{n+1})^*$ \ g\'en\'erique on aura \ $\{ l = 0 \} \cap S = \{0\} $. Si, par exemple, \ $(\frac{\partial f}{\partial x_1}, \dots, \frac{\partial f}{\partial x_n})$ \ est une suite r\'eguli\`ere en \ $0$, toute forme lin\'eaire d'un ouvert dense dans un voisinage ouvert assez petit  de la forme lin\'eaire \ $ x \mapsto x_0$ \ sera telle que le lieu critique de l'application \ $(f,l)$ \ sera de dimension 1 et n'aura pas de composante irr\'eductible contenue dans \ $\{l = 0\} $. On aura donc un ouvert dense sur lequel les deux premi\`eres conditions de l'\'enonc\'e sont r\'ealis\'ees avec, de plus, le fait que chaque fonction \ $f_t$ \ ne pr\'esente que des points singuliers isol\'es dans l'hyperplan  \ $\{l = t\}$ \ pour tout \ $t$ \ assez voisin de \ $0$.\\
Le fait que \ $\mu$ \ reste localement  constant sur \ $S^*$ \ au voisinage de l'origine est alors cons\'equence du fait que le faisceau  \ $l_*(\Omega_{/l}^n\big/ d_{/l}f \wedge \Omega_{/l}^{n-1})$ \ est coh\'erent sur \ $\mathcal{O}_D$ \ et donc localement libre en dehors de l'origine, puisque la restriction de \ $l$ \ au support de ce faisceau coh\'erent est finie. $\hfill \blacksquare$

\bigskip

\subsection{Rappels.}

Rappelons maintenant quelques constructions et quelques r\'esultats donn\'es dans [B.II]  pour l'\'etude de la "situation consid\'er\'ee".

\parag{R.1} D\'efinissons le faisceau des formes $t-$relatives
$$ \Omega^{\bullet}_/ : = \Omega^{\bullet}_X \big/ dt \wedge \Omega^{\bullet -1}_X $$
la diff\'erentielle $t-$relative \'etant induite par la diff\'erentielle de de Rham sur le quotient. Comme la fonction \ $t$ \ est non singuli\`ere on a un scindage naturel
$$ \Omega^{\bullet}_X \simeq dt\wedge \Omega^{\bullet-1}_X \oplus \Omega^{\bullet}_/ $$
ainsi qu'un isomorphisme diff\'erentiel gradu\'e \ 
$$\wedge dt : \Omega^{\bullet}_/ \to dt\wedge \Omega^{\bullet}_X .$$
Il est int\'eressant de remarquer que dans la suite exacte de cohomologie de la suite exacte courte de complexes
$$ 0 \to (\Omega_/^{\bullet}, d_/^{\bullet})[-1] \overset{\wedge dt}{\to} (\Omega_X^{\bullet}, d^{\bullet}) \to (\Omega_/^{\bullet}, d_/^{\bullet}) \to 0 $$
le connecteur s'identifie \`a la d\'erivation \ $\partial / \partial t $.

\smallskip

On a exactitude en degr\'es strictement positifs du complexe \ $(\Omega^{\bullet}_/, d^{\bullet}_/)$.\\
On a \'egalement exactitude en degr\'es \ $d \in [1,n-1]$ \ du complexe \ $(\Omega^{\bullet}_/, (\wedge d_/f)^{\bullet})$ \ gr\^ace \`a la finitude de \ $ t : \Sigma \to D$.

\smallskip

Ces propri\'et\'es  s'etendent imm\'ediatement au compl\'et\'e formel \ $\hat{\Omega}^{\bullet}$ \ en \ $f$.

\parag{R.2} Introduisons maintenant le sous-complexe \ $(\hat{K}er\,d_/f)^{\bullet}, d^{\bullet}_/)$ \ du complexe \ $(\hat{\Omega}^{\bullet}_/, d^{\bullet}_/)$ \ o\`u \ $(\hat{K}er\,d_/f)^{\bullet}$ \ est le noyau de la multiplication ext\'erieure  
$$\wedge d_/f : \hat{\Omega}^{\bullet}_/ \to \hat{\Omega}^{\bullet+1}_/ .$$
Nous noterons respectivement par \ $\mathcal{E}$ \ et \ $\mathbb{E}$ \ les faisceaux de cohomologie de ce complexe en degr\'es $1$ \ et \ $n$. Ils sont naturellement munis d'op\'erations \ $a$ \ et $b$ \ v\'erifiant \ $a.b - b.a = b^2$ \ d\'eduites de la multiplication par \ $f$ \ et de \ $\wedge d_/f \circ (d_/)^{-1}$. Celles-ci commutent \`a l'action naturelle de \ $t^{-1}(\mathcal{O}_D)$ \ sur ces faisceaux.\\
On remarquera que la compl\'etion formelle fait que ces faisceaux, \`a priori port\'es par \ $\Sigma$, sont \`a support dans \ $S$ \ puisque \ $S = \Sigma \cap \{f = 0 \}$, gr\^ace \`a  notre hypoth\`ese sur \ $t$. En effet un germe de courbe irr\'eductible \ $(\Gamma,0)$ \ contenu dans \ $\Sigma \cap \{ f = 0 \}$ \ v\'erifie ou bien \ $(\Gamma,0) \subset (S,0)$ \ ou bien \ $\Gamma \cap S = \{0\}$. Dans ce dernier cas, \ $\Gamma \setminus \{0\}$ \ est connexe et contenue dans l'ouvert lisse \ $\{ df \not= 0 \}$ \ de \ $\{ f = 0 \}$, et la fonction \ $t_{\vert \Gamma \setminus \{0\}}$ \ est localement constante, donc identiquement nulle. On en d\'eduit que \ $(\Gamma,0) \subset (\Sigma,0) \cap \{ t = 0 \} = \{0\}$, ce qui est absurde. D'o\`u notre assertion.\\

\smallskip

Le th\'eor\`eme suivant a \'et\'e \'etabli dans [BII].

\begin{thm}( [B.II] th. 4.2.1.)
Dans la "situation consid\'er\'ee" au d\'ebut de cette section, les faisceaux de cohomologie du complexe \ $(\hat{K}er\,d_/f)^{\bullet}, d^{\bullet}_/)$ \ sont nuls en degr\'es diff\'erents de \ $1$ \ et \ $n$. En degr\'e $1$ \ et \ $n$ \ les faisceaux \ $\mathcal{E}$ \ et \ $\mathbb{E}$ \ sont  des \ $t^{-1}(\mathcal{O}_D[[b]])-$modules coh\'erents, localement libres sur \ $S^*$, si le disque \ $D$ \ est assez petit.
\end{thm}

\parag{R.3} Dans la "situation consid\'er\'ee"  d\'efinissons le morphisme de faisceaux de \ $\mathbb{C}-$espaces vectoriels sur \ $S$
$$ \nabla :  \mathbb{E} \to \mathbb{E}$$
en posant 
 $$\nabla[d_/\xi] : = [ d_/f \wedge \frac{\partial \xi}{\partial t } - \frac{\partial f}{\partial t }.d\xi ].$$
 La v\'erification du fait que ce morphisme de faisceaux est bien d\'efini est facile. Les propri\'et\'es suivantes sont d\'emontr\'ees dans [B.II] lemme 4.2.5 et  proposition 4.2.8.
 
 \begin{prop}
 Soit \ $\mathbb{P} : = \{ x \in  \mathbb{E} \ / \ \nabla(x) \in b.\mathbb{E} \}$. C'est un sous-$t^{-1}(\mathcal{O}_D[[b]])-$module de \ $\mathbb{E}$ \ qui est stable par \ $a$. Le morphisme de faisceaux
 $$ b^{-1}.\nabla : \mathbb{P} \to \mathbb{E} $$
 induit par \ $\nabla$ \ est une \ $t^{-1}(\mathcal{O}_D)-$connexion qui commute \`a \ $a$ \ et \ $b$.\\
 Son noyau et son conoyau sont respectivement canoniquement isomorphes comme \ $\hat{\mathcal{A}}-$modules aux faisceaux de cohomologie  \ $\hat{\mathcal{H}}^n$ \ et \ $\hat{\mathcal{H}}^{n+1}$ \ du complexe \ $(\hat{K}er\,df^{\bullet}, d^{\bullet})$.
 \end{prop}
 
 \bigskip
 
L'alg\`ebre \ $\hat{\mathcal{A}}$ \ qui n'intervient ici essentiellement que pour traduire la compatibilit\'e aux op\'erations \ $a$ \ et \ $b$ \ de cet isomorphisme, a \'et\'e d\'efinie (bri\`evement) au 1.2.1. Pour plus de d\'etails  sur cette alg\`ebre nous renvoyons le lecteur \`a [B.95].
 
 \subsection{La \ $b^{-1}.\nabla-$finitude.}
 
 Le contr\^ole du noyau de \ $\nabla$ \ sera obtenu gr\^ace \`a la propri\'et\'e suivante.
 
   \begin{defn}\label{Nabla-fin.}
   On dira que \ $\mathbb{E}$ \ est {\bf \ $b^{-1}\nabla-$fini} si, localement sur \ $D^*$, il existe un entier \ $N$ \ et  un sous-$t^{-1}(\mathcal{O}_D[[b]])-$module coh\'erent  \ $\hat{G}$ \ de \ $\mathbb{P}$ \ qui est stable par \ $b^{-1}\nabla$ \ et contient \ $b^N.\mathbb{E}$.
   \end{defn}
   
   \begin{prop}\label{G}
  Il existe sur \ $D$ \ un plus grand sous-$t^{-1}(\mathcal{O}_D[[b]])-$module \ $\mathbb{G}$ \ de \ $\mathbb{P}$, stable par  \ $b^{-1}\nabla$.\\
   Il v\'erifie de plus les propri\'et\'es suivantes :
   \begin{enumerate}[(1)]
   \item On a \ $\mathbb{G} = \{ x \in \mathbb{E} \ / \ \forall \nu \in \mathbb{N} \quad \nabla^{\nu}(x) \in b^{\nu}.\mathbb{E}\}.$
     \item Si \ $x \in \mathbb{E}$ \ v\'erifie \ $b^{-1}\nabla(x) \in \mathbb{G} + t^{-1}(\mathcal{O}_D[[b]]).x$ \ alors \ $x \in  \mathbb{G}$.\\
      En particulier on a  \ $Ker\nabla \subset \mathbb{G}$.
    \item Si des germes \ $\varphi \in \mathcal{O}_D \setminus \{0\}$ \ et \ $x \in \mathbb{E}$ \ sont tels que \ $t^{-1}(\varphi).x \in \mathbb{G}$ \ alors \ $x \in \mathbb{G}$.
     \item Si des germes \ $\varphi \in \mathcal{O}_D \setminus \{0\}$ \ et \ $x \in \mathbb{G}$ \ sont tels que \ $t^{-1}(\varphi).x \in b.\mathbb{G}$ \ alors \ $x \in b.\mathbb{G}$.
   \end{enumerate}
     Supposons maintenant que \ $\mathbb{E} $ soit \ $b^{-1}\nabla-$fini sur \ $S^*$. Alors \ $\mathbb{G}$ \ est \ $t^{-1}(\mathcal{O}_D[[b]])$ \ localement libre de rang fini sur \ $S^*$ et il est stable par \ $a$. De plus il v\'erifie les deux propri\'et\'es suivantes:
     \begin{enumerate}[i)]
       \item Le sous-faisceau \ $Ker\, \nabla$ \ est localement constant sur \ $S^*$, de fibre un (a,b)-module r\'egulier g\'eom\'etrique.
    \item Soit \ $U \subset D^*$ \  un ouvert connexe et simplement connexe, et soit  \ $V \subset S^*$ \ une composante connexe de \ $t^{-1}(U)$. On a un isomorphisme 
   $$ \Gamma(V, \mathbb{G}) \simeq \Gamma(V, Ker\, \nabla/b.Ker\, \nabla) \otimes_{\mathbb{C}} t^{-1}(\mathcal{O}_{D^*}(U)[[b]]).$$
   \end{enumerate}
 \end{prop}
 
 \parag{Remarque} Pr\'ecisons que \ $\mathbb{G}$ \ est maximal au sens suivant : 
 pour tout ouvert \ $\Omega \subset D$ \ et tout sous-faisceau \ $\Gamma $ \ de \ $\mathcal{O}_D[[b]]-$modules de \ $\mathbb{P}_{\vert \Omega}$ \ stable par \ $b^{-1}.\nabla$ \ on aura  \ $\Gamma \subset \mathbb{G}$. $\hfill \square$

   \parag{Preuve} Notons \ $\mathbb{G}_k : = \{ x \in \mathbb{E} \ / \ \forall \nu \in [0,k] \quad \nabla^{\nu}(x) \in b^{\nu}.\mathbb{E}\}$, pour \ $k \in \mathbb{N}$. On a donc \ $\mathbb{G} = \cap_{k \geq 0} \ \mathbb{G}_k$. Comme \ $\mathbb{G}_k$ \ est manifestement stable par \ $b$, pour voir que \ $\mathbb{G}_k$ \ est un sous-$\mathcal{O}_D[[b]]-$module, il suffit de voir qu'il est stable par l'action de \ $t^{-1}(\mathcal{O}_D)$.\\
 Pour \ $\varphi \in \mathcal{O}_D$ \ et \ $x \in \mathbb{E}$ \ on a (en oubliant de noter l'image r\'eciproque par \ $t$)
  \begin{equation*}
   \nabla^{\nu}(\varphi.x) = \sum_{j=0}^{\nu} \ C_{\nu}^j.\varphi^{(j)}.b^j.\nabla^{\nu - j}(x) \tag{@}
   \end{equation*}
    o\`u \ $\varphi^{(j)}$ \ d\'esigne la d\'eriv\'ee \ $j-$\`eme de \ $\varphi$.\\
  Cette formule montre que si  \ $x \in \mathbb{G}_k$ \ on aura \ $\varphi.x \in \mathbb{G}_k$.\\
   L'inclusion de tout sous-$t^{-1}(\mathcal{O}_D[[b]])-$module stable par \ $b^{-1}\nabla$ \  dans \
    $\mathbb{G}$ \  est imm\'ediate. Pour prouver (1)  il suffit donc de montrer que \ $\mathbb{G}$ \  est stable par \ $b^{-1}\nabla$. Or 
    $$\nabla(\mathbb{G}_k) = \{ \nabla(x) \ / x \in \mathbb{G}_k\}\subset \{y \in \mathbb{E} \ / \ \nabla^{\nu}(y) \in b^{\nu+1}\mathbb{E} \quad \forall \nu \in [0,k-1]\} $$
      et ce dernier groupe est \'egal \`a 
       $$  \{bz \ / \ z \in \mathbb{E} \ / \ \nabla^{\nu}(z) \in b^{\nu}.\mathbb{E}\quad \forall \nu \in [0,k-1]\} = b.\mathbb{G}_{k-1}.$$
    On a donc bien \ $\nabla(\mathbb{G}) \subset b.\mathbb{G}$.
    
    \smallskip
    
    Si \ $x \in \mathbb{E}$ \ v\'erifie \ $b^{-1}\nabla(x) \in \mathbb{G} + t^{-1}(\mathcal{O}_D[[b]]).x$ \ consid\'erons   \ $G : = t^{-1}(\mathcal{O}_D[[b]]).x + \mathbb{G}$. On a alors, pour \ $y \in \mathbb{G}$, 
   $$ \nabla(t^{-1}(\alpha).x + y) \in  t^{-1}(\alpha').bx +  t^{-1}(\alpha).\nabla(x) + b.\mathbb{G} \subset b.(t^{-1}(\mathcal{O}_D[[b]]).x + \mathbb{G}) $$
   o\`u \ $\alpha'$ \ d\'esigne la d\'eriv\'ee en \ $t$ \ de \ $\alpha \in \mathcal{O}_D[[b]]$. Ceci  montre que \ $G = \mathbb{G}$ \ et donc que \ $x \in \mathbb{G}$. Ceci montre (2).\\
   Si maintenant on a \ $t^{-1}(\varphi).x \in \mathbb{G}$ \ avec \ $\varphi \in \mathcal{O}_D \setminus \{0\}$ \ et \ $x \in \mathbb{E}$ \ montrons par r\'ecurrence sur \ $\nu \in  \mathbb{N}$ \ que l'on a \ $\nabla^{\nu}(x) \in b^{\nu}.\mathbb{E}$. La propri\'et\'e \'etant claire pour \ $\nu = 0$ \ supposons-la vraie pour \ $\nu -1$ \ et montrons-la pour \ $\nu$. Comme on a la relation \ $(@)$
    on d\'eduit de l'hypoth\`ese \ $t^{-1}(\varphi).x \in \mathbb{G}$ \ et de l'hypoth\`ese de r\'ecurrence que \ $t^{-1}(\varphi).\nabla^{\nu}(x) \in b^{\nu}.\mathbb{E}$. Mais sur \ $D^*$ \ le faisceau \ $\mathbb{E}/b^{\nu}.\mathbb{E}$ \ est localement \ $t^{-1}(\mathcal{O}_D)-$libre\footnote{car c'est le cas pour \ $\nu = 1$ \ et donc pour tout \ $\nu$ \ en raisonnant par r\'ecurrence sur la suite exacte 
    $$ 0 \to b^k.\mathbb{E}\big/b^{k+1}.\mathbb{E} \to \mathbb{E}\big/b^{k+1}.\mathbb{E} \to \mathbb{E}\big/b^{k}.\mathbb{E} \to 0. $$
    puisque \ $b^k : \mathbb{E}\big/ b.\mathbb{E} \to b^k.\mathbb{E}\big/b^{k+1}.\mathbb{E}$ \ est un isomorphisme de faisceaux de \ $\mathcal{O}_D-$modules.} on en d\'eduit que \ $\nabla^{\nu}(x) \in b^{\nu}.\mathbb{E}$. On en conclut que \ $x \in \mathbb{G}$, ce qui ach\`eve la preuve du point (3).\\
    Le point (4) se montre de la m\^eme fa{\c c}on en utilisant l'\'egalit\'e (d\'ej\`a vue plus haut)
    $$ b.\mathbb{G} = \{ x \in \mathbb{E} \ / \ \forall \nu \in \mathbb{N} \quad \nabla^{\nu}(x) \in b^{\nu+1}.\mathbb{E}\}.$$

    \smallskip
    
    Supposons maintenant que \ $\mathbb{E}$ \ est \ $b^{-1}\nabla-$fini. \\
    Comme nos assertions sont locales sur \ $S^*$, nous pouvons nous placer au-dessus d'un disque \ $U \subset D^*$,  consid\'erer une composante connexe \ $V$ \ de \ $t^{-1}(U)$ \ et supposer que l'on ait sur \ $V$ \ un sous-$\mathcal{O}_D[[b]]-$module coh\'erent \ $\hat{G}$ \ contenant \ $b^N.\mathbb{E}_{\vert V}$ \ et stable par \ $b^{-1}.\nabla$.\\
    Pour simplifier les notations, nous identifierons \ $U$ \ et \ $V$ \ via \ $t$.
    
      \parag{\'Etape 1}On se ram\`ene au cas o\`u  le \ $\mathcal{O}_D-$module \ $\mathbb{P}/\hat{G}$ \ est libre   sur \ $V$ \ et de rang  not\'e \ $p$,  et o\`u, de plus, il est facteur direct du  \ $\mathcal{O}_D-$module libre \ $\mathbb{E}/\hat{G}$ \ dont le  rang sur \ $V$ \ sera not\'e \ $p + q$. Ceci est montr\'e dans le lemme suivant.
    
 \begin{lemma}
    On se place dans la situation pr\'ecis\'ee ci-dessus. Il existe un sous-$\mathcal{O}_D[[b]]-$module \ $\tilde{G}$ \ coh\'erent stable par \ $b^{-1}.\nabla$ \ qui contient \ $\hat{G}$ \ et qui v\'erifie de plus les propri\'et\'es suivantes :
    \begin{enumerate}[a)]
    \item  Le \ $\mathcal{O}_D-$module \ $\mathbb{P}/\tilde{G}$ \ est libre.
    \item  La suite exacte de \ $\mathcal{O}_D-$modules
    $$ 0 \to \mathbb{P}/\tilde{G} \to \mathbb{E}/\tilde{G} \to \mathbb{E}/\mathbb{P} \to 0 $$
    est scind\'ee.
    \end{enumerate}
    \end{lemma}
    
    \parag{Preuve} Le quotient \ $\mathbb{P}/\hat{G}$ \ est \ $\mathcal{O}_D-$coh\'erent puisque \ $\hat{G}$ \ contient \ $b^N.\mathbb{P}$. Donc sa \ $\mathcal{O}_D-$torsion \ $T$ \  est coh\'erente sur \ $\mathcal{O}_D$. Si \ $\pi : \mathbb{P} \to \mathbb{P}/\hat{G}$ \ est l'application quotient, d\'efinissons
    $$ \tilde{G} : \pi^{-1}(T).$$
    Il est \ $\mathcal{O}_D[[b]]-$coh\'erent comme noyau et la propri\'et\'e a) est \'evidente. Il contient \ $\hat{G}$ \ par d\'efinition. Montrons la stabilit\'e par \ $b^{-1}.\nabla$. \\
    Le probl\`eme est local au voisinage de chaque point \ $t_0$ \ appartenant au support du faisceau \ $T$, support qui est ferm\'e et discret dans \ $V$. Pr\`es d'un tel point, on aura \ $x \in \tilde{G}_{t_0}$ \ si et seulement si  il existe un entier \ $k \geq 1$ \ tel que l'on ait \  $(t-t_0)^k.x \in \hat{G}$. Alors l'\'egalit\'e 
     $$ (t-t_0)^{k+1}.b^{-1}.\nabla(x) = b^{-1}.\nabla((t-t_0)^{k+1}.x) - (k+1).(t-t_0)^k.x$$
     \ montre que \ $b^{-1}.\nabla(x) \in \tilde{G}$ \ puisque \ $\hat{G}$ \ est stable par \ $b^{-1}.\nabla$.\\
     Le scindage de la suite exacte est cons\'equence du fait que le \ $\mathcal{O}_D-$module \ $\mathbb{E}/\mathbb{P}$ \ est libre. En effet si on a \ $\varphi.x \in \mathbb{P}$ \ avec \ $\varphi \not= 0 $, cela signifie que \ $\nabla(\varphi.x) = \varphi'.bx + \varphi.\nabla(x) \in b.\mathbb{E}$, c'est \`a dire que \ $\varphi.\nabla(x) \in b.\mathbb{E}$. Comme \ $\mathbb{E}/b.\mathbb{E}$ \ est sans \ $\mathcal{O}_D-$torsion, cela signifie que \ $\nabla(x) \in b. \mathbb{E}$ \ c'est \`a dire que \ $x \in \mathbb{P}$. $\hfill \blacksquare$
     
     \bigskip
     
     On remplacera dans la suite \ $\hat{G}$ \ par \ $\tilde{G}$ \ mais en continuant \`a le noter \ $\hat{G}$.
     
     \smallskip

   Ayant choisi ces trivialisations, la connexion
    $$ b^{-1}.\nabla : \mathbb{P}/\hat{G} \to \mathbb{E}/\hat{G}$$
     se lit comme une \ $\mathcal{O}_D-$connexion associ\'ee \`a l'inclusion naturelle de \ $\mathcal{O}^p_D $ \ dans \ $\mathcal{O}^p_D \oplus \mathcal{O}^q_D$ 
    $$ \partial : \mathcal{O}^p_D \to \mathcal{O}^p_D \oplus \mathcal{O}^q_D $$
    sur le disque \ $V$.
    
    \parag{Deuxi\`eme \'etape} Elle est donn\'ee par le lemme simple suivant.
    
    \begin{lemma}\label{Ker}
    Soit
    $$ \partial : \mathcal{O}^p_D \to \mathcal{O}^p_D \oplus \mathcal{O}^q_D $$
une \ $\mathcal{O}_D-$connexion sur un disque \ $V \subset D$. Alors \ $Ker \partial$ \ est un sous-faisceau constant de \ $\mathbb{C}_D-$modules de \ $\mathcal{O}_D^p$ \ de rang au plus \'egal \`a \ $p$ \ et le plus grand sous-faisceau de \ $\mathcal{O}_D-$module stable par \ $\partial$ \ est \'egal \`a \ $\mathcal{O}_D.Ker \partial  $. Il est donc libre et facteur direct dans \ $\mathcal{O}^p_D$.
\end{lemma}

\parag{Preuve} Soit \ $\pi :  \mathcal{O}^p_D \oplus \mathcal{O}^q_D \to \mathcal{O}^p_D$ \ la projection. Alors \ $\pi \circ \partial$ \ est une connexion sur \ $\mathcal{O}^p_D$ \ au-dessus du disque \ $V$. Soit \ $e : = (e_1, \cdots, e_p)$ \ une base horizontale pour cette connexion. Notons \ $\varepsilon : = (\varepsilon_1, \cdots, \varepsilon_q)$ \ la base canonique de \ $\mathcal{O}^q_D $. D\'efinissons alors la matrice holomorphe \ $M : V \to L(\mathbb{C}^p, \mathbb{C}^q)$ \ par la formule
$$ \partial e = M.\varepsilon .$$
Pour \ $X \in \mathcal{O}_D^p $ \ on aura
$$ \partial(^tX.e) = \ ^tX'.e \oplus \ ^tX.M.\varepsilon \ ,$$
o\`u \ $X'$ \ d\'esigne la d\'eriv\'ee de \ $X$. On en d\'eduit que \ $^tX.e \in Ker \partial$ \ si et seulement si \ $X \in \mathbb{C}^p$ \ et v\'erifie \ $^tX.M \equiv 0 $ \ sur \ $V$. Ceci d\'efinit un sous-espace vectoriel de \ $\mathbb{C}^p$ \ et montre donc la premi\`ere assertion.\\
Les autres assertions sont imm\'ediates. $\hfill \blacksquare$

\parag{Troisi\`eme \'etape} Posons maintenant, avec les notations pr\'ec\'edentes
$$ G : = \mathcal{O}_D[[b]].Ker \partial  + \hat{G} .$$
Alors \ $G$ \ est \ $\mathcal{O}_D[[b]]-$coh\'erent comme somme de deux sous-faisceaux coh\'erents de \ $\mathbb{P}$ \ qui est coh\'erent (voir [B.II] 7.2).
Comme \ $\mathbb{G}/\hat{G}$ \ est un sous-$\mathcal{O}_D-$module de \ $\mathbb{P}/\hat{G}$ \ qui est stable par \ $\partial$, on aura, d'apr\`es le lemme \ref{Ker}
$$ \mathbb{G} \subset G .$$
Mais \ $G$ \ est un sous-$\mathcal{O}_D[[b]]-$module de \ $\mathbb{P}$ \ stable par \ $b^{-1}.\nabla$. En effet, si \ $x \in G$, comme on a
$$ b^{-1}.\nabla(x) \in b^{-1}.\nabla(\hat{G}) + \mathcal{O}_D.Ker \partial + \hat{G}  \subset G.$$
On en conclut  que \ $G = \mathbb{G}$ \ sur l'ouvert \ $V$. Ceci donne alors la coh\'erence de \ $\mathbb{G}$ \ sur \ $\mathcal{O}_D[[b]]$. Comme le faisceau \ $\mathbb{G}/b.\mathbb{G}$ \ est sans \ $\mathcal{O}_D-$torsion d'apr\`es la propri\'et\'e (4) prouv\'ee plus haut, on en conclut que \ $\mathbb{G}$ \ est (localement) libre de rang fini sur \ $\mathcal{O}_D[[b]]$. De plus, d'apr\`es la propri\'et\'e (2) il contient \ $Ker \nabla$.\\
Pour conclure, il suffit alors d'appliquer le th\'eor\`eme de Cauchy \`a la \ $\mathcal{O}_D[[b]]-$connexion \ $b^{-1}.\nabla$ \ de \ $\mathbb{G}$. $\hfill \blacksquare$

 \section{Int\'egration "\`a la Malgrange".}
 
 Ce chapitre est, bien s\^ur, directement inspir\'e de [M. 74].
     
     \subsection{Complexe de de Rham absolu et \ $b^{-1}.\nabla$.}
     
     Le lien entre le complexe de de Rham absolu \ $(\hat{K}er\ df)^{\bullet}, d^{\bullet})$ \ et la connexion \ $b^{-1}.\nabla$ \ introduite plus haut est pr\'ecis\'e par le lemme suivant     
       \begin{lemma}
    Soit \ $n$ \ un entier au moins \'egal \`a \ $ 2$. On a un  morphisme de complexes
   $$ \xymatrix{ \cdots \ar[r] & (\hat{K}er\ df)^{n-1} \ar[r]^d \ar[d] & (\hat{K}er\ df)^n \ar[r]^d \ar[d]^{\alpha} & (\hat{K}er\ df)^{n+1} \ar[r] \ar[d]^{\beta} & 0 \\
                            \cdots \ar[r] &  0 \ar[r]                               & \mathbb{P} \ar[r]^{b^{-1}\nabla}             & \mathbb{E}   \ar[r]        & 0}  $$
 o\`u \ $\alpha$ \ est induite par la projection "\'evidente" de \ $\hat{\Omega}^n$ \ sur \ $\hat{\Omega}^n_/$ \ et o\`u \ $\beta$ \ est l'inverse de l'isomorphisme \ $ \wedge dt : \hat{\Omega}^n_/ \to \hat{\Omega}^{n+1}$. Il  induit des isomorphismes (a,b)-lin\'eaires sur les faisceaux de cohomologie, si on remplace \ $(\hat{K}er\ df)^0$ \ par \ $(\hat{K}er\ df)^1 \cap Ker\ d $ \ en degr\'e 0 dans le premier complexe (avec l'inclusion \'evidente).
   \end{lemma}
   
   \bigskip
   
   \noindent \textit{Preuve.} V\'erifions d\'ej\`a les \'egalit\'es \ $\alpha \circ d = 0 $ \ et \ $\beta \circ d = b^{-1}\nabla \circ \alpha$.\\
    Si \ $u + dt\wedge v$ \ est dans \ $(\hat{K}er\ df)^{n-1}$ \ avec \ $u \in \hat{\Omega}^{n-1}_/ $ \ et \ $v \in \hat{\Omega}^{n-2}_/ $, on aura \ $d_/f \wedge u = 0$ \ et   \ $(\alpha\circ d)(u + dt\wedge v) = d_/u$. On trouve donc bien z\'ero dans \ $\mathbb{P} \subset \mathbb{E}$.\\
   Pour \ $d_/\xi + dt\wedge \eta \in (\hat{K}er\ df)^n$ \ on aura 
   $$(\beta\circ d)(d_/\xi + dt\wedge \eta) = \beta(dt \wedge( \frac{\partial d_/\xi}{\partial t} - d_/\eta)) = \frac{\partial d_/\xi}{\partial t} - d_/\eta =d_/(\frac{\partial \xi}{\partial t} - \eta).$$
   La condition \ $df \wedge (d_/\xi + dt\wedge \eta) = 0 $ \ donne \ $ \frac{\partial f}{\partial t}.d_/\xi = d_/f \wedge \eta $ \ et on a donc \ 
   $$\nabla(\alpha(d_/\xi + dt\wedge \eta)) = \nabla(d_/\xi) = d_/f \wedge (\frac{\partial \xi}{\partial t} - \eta) $$
   ce qui donne bien l'\'egalit\'e d\'esir\'ee.\\
   La fin de la preuve consiste alors \`a v\'erifier que les morphismes induits en cohomologie co{\"i}ncident avec ceux de la proposition \ref{connexion}. Cette v\'erification simple est laiss\'ee au lecteur. $\hfill \blacksquare$

   \subsection{D\'erivations  d'int\'egrales}
   
     Comme pr\`es du point \ $p \in S^*$ \ la famille des fonctions \`a singularit\'es isol\'ees donn\'ee par \ $t \to f_t$ \ est \`a \ $\mu-$constant elle est topologiquement triviale au voisinage de \ $p$. En identifiant \ $S^*$ \ et \ $D^*$ \ via la fonction \ $t$ \ au voisinage de \ $p$, on a donc l'existence d'un voisinage \ $U$ \ de \ $p$ \ dans \ $X$ et une application continue \ $\Phi : U \to U_p : = t^{-1}(t(p)) \cap U $ \ donnant le diagramme  commutatif
  $$ \xymatrix{ U \ar[d]^{f\times t} \ar[r]^{\Phi\times t} & U_p \times \Delta \ar[d]^{f\times Id}\\
  D \times \Delta \ar[r]^{=} & D \times \Delta} $$
  o\`u \ $(\Phi\times t)$ \ est un hom\'eomorphisme et \ $U_p$ \ une boule de Milnor en \ $p$ \ pour la restriction \ $f_p$ \ de \ $f$ \ \`a l'hyperplan  \ $ t^{-1}(t(p))$.\\
 Soit \ $\gamma \in H_{n-1}(\{f = \varepsilon\} \cap U_p, \mathbb{C})$, pour \ $\varepsilon > 0$ \ assez petit. On a pour \ $0 < \vert s \vert < \varepsilon$ \ et \ $t$ \ assez voisin de \ $t(p))$ \ une famille horizontale multiforme de \ $(n-1)-$cycles compacts  contenus dans les fibres \ $ \{f(t,x) = s\} \cap U$.\\
  \'Etant donn\'ee une section \ $\omega$ \ du faisceau \ $\mathbb{E}$ \ sur un voisinage ouvert de \ $p \in S^*$, nous d\'efinirons les int\'egrales "\`a la Malgrange"  
  $$ \int_{\gamma_{s,t}} \ \frac{\omega}{d_/f} .$$
  On remarquera d\'ej\`a qu'une telle int\'egrale ne d\'epend que de la classe de \ $\omega$ \ dans \ $\mathbb{E}$ \ en raison de la formule de Stokes, puisque 
   $$\mathbb{E} = \hat{\Omega}^n_/\big/ d_/(\hat{K}er\,d_/f^{n-1}) =  \hat{\Omega}^n_/\big/ d_/f\wedge d_/\hat{\Omega}^{n-2}_/ .$$
  La proposition suivante exprime les d\'eriv\'ees partielles en \ $s$ \ et \ $t$ \ de ces fonctions holomorphes en \ $(s,t)$ \ qui sont multivalu\'ees en $s$.

      \begin{prop}\label{formules}
   Pour tout \ $d_/\xi \in \mathbb{E} $ \ on a la formule de d\'erivation "en $s$"
   $$ \frac{\partial}{\partial s} \big( \int_{\gamma_{s,t}} \ \xi \big) = \int_{\gamma_{s,t}} \frac{d_/\xi}{d_/f}  $$
   que l'on peut lire sous la forme
   $$primitive \ "en \ s" \  de \ \Big(\int_{\gamma_{s,t}} \ \frac{d_/\xi}{d_/f} \Big)=  \int_{\gamma_{s,t}} \ \frac{d_/f \wedge \xi}{d_/f} =  \int_{\gamma_{s,t}} \  \frac{b(d_/\xi)}{d_/f} \ .$$
   Pour tout \ $d_/\xi \in \mathbb{E}$ \ on a la formule de d\'erivation "en $t$"
   $$\frac{\partial}{\partial t} \big( \int_{\gamma_{s,t}}\ \xi \big) = \int_{\gamma_{s,t}} \ \frac{\nabla(d_/\xi)}{d_/f} $$
   que l'on peut lire sous la forme
   $$\frac{\partial}{\partial t}\Big(\int_{\gamma_{s,t}} \ \frac{d_/\xi}{d_/f} \Big)=  \int_{\gamma_{s,t}} \ \frac{b^{-1}\nabla(d_/\xi)}{d_/f}\ $$
   pour \ $d_/\xi \in \mathbb{P}$.

   \end{prop}
   
   \bigskip
   
   \noindent \textit{Preuve.} Soit \ $\Delta \subset \mathbb{H}\times \mathbb{C}$ \ un polydisque. L'hypoth\`ese de trivialit\'e topologique permet de trouver \ $\psi \in \mathcal{C}^{\infty, n-1}_{c/f}(f^{-1}(\Delta))$ \ une \ $(n-1)-$forme $d-$ferm\'ee induisant pour chaque \ $(s,t) \in \Delta$ la classe d\'efinie par \ $\gamma_{s,t}$ \ dans 
   $$H^{n-1}_c(f^{-1}_t(s)\cap U, \mathbb{C}) \simeq H_{n-1}(f_t^{-1}(s)\cap U, \mathbb{C}).$$
   
    Pour \ $\xi \in \Omega^{n-1}_/ $ \ posons
   $$ F(s,t) : = \int_{\gamma_{s,t}} \ \xi \  = \int_{f_t^{-1}(s)} \ \xi \wedge \psi \quad \quad \forall (s,t) \in \Delta .$$
   Commen{\c c}ons par prouver l'holomorphie de \ $F$ \ sur \ $\Delta$. Comme \ $F$ \ est manifestement une fonction continue, il suffit de prouver que son \ $\bar{\partial}$ \ au sens des distributions est nul. Consid\'erons alors une forme test \ $\theta \in \mathcal{C}^{\infty, (2,1)}_c(\Delta) $. On a
   $$ < \bar{\partial}F, \theta > = - < F, \bar{\partial}\theta > = -\int_{f^{-1}(\Delta)} \  \xi \wedge \psi \wedge d(f^*(\theta)) $$
   puisque \ $ \bar{\partial}\theta = d\theta$ \ et puisque loin du lieu critique de \ $f$ \ le th\'eor\`eme de Fubini banal s'applique \`a des fonction continues. Comme on a 
   $$ d(\xi \wedge \psi \wedge f^*(\theta)) = \xi \wedge \psi \wedge d(f^*(\theta))) $$
   puisque \ $d\psi = 0 $ \ et que \ $d\xi \wedge f^*(\theta)$ \ est de type \ $(n+2,1)$ \ donc nulle, la formule de Stokes permet de conclure \`a l'holomorphie de \ $F$.\\
   Pour calculer \ $\frac{\partial F}{\partial s}$ \ au sens des distributions, consid\'erons maintenant une forme test \ $\zeta \in \mathcal{C}^{\infty, (0,2)}_c(\Delta)$. On aura
   $$ < \frac{\partial F}{\partial s}.ds , dt \wedge\zeta > = - < F , d(dt \wedge \zeta)> =  - \int_{f^{-1}(\Delta)} \ \xi \wedge \psi \wedge d(f^*(dt \wedge \zeta )) .$$
   Comme on a \ $d\xi = d_/\xi + dt \wedge \frac{\partial \xi}{\partial t} $ \ la formule de Stokes donne
   $$ < \frac{\partial F}{\partial s}.ds , dt \wedge\zeta > = \int_{f^{-1}(\Delta)} \ d_/ \xi \wedge \psi \wedge f^*(dt \wedge\zeta) .$$
   Le th\'eor\`eme de Fubini donne alors
   $$  < \frac{\partial F}{\partial s}.ds , dt \wedge\zeta> =  \int_{\Delta} \big(\int_{f_t^{-1}(s)} \ \frac{d_/\xi}{d_/f} \wedge \psi \big). ds \wedge dt \wedge \zeta $$
   ce qui donne bien notre formule de d\'erivation "en s".\\
 Soit \ $\zeta$ \ une forme-test choisie comme plus haut. On a
   $$  < \frac{\partial F}{\partial t}.dt ,  ds \wedge \zeta > = - < F, d(ds \wedge \zeta) > = - \int_{f^{-1}(\Delta)} \ \xi \wedge \psi \wedge d(f^*(ds \wedge \zeta)) .$$
   La formule de Stokes donne alors
   $$ < \frac{\partial F}{\partial t}.dt ,  ds \wedge \zeta > =  \int_{f^{-1}(\Delta)} \ \frac{\partial \xi}{\partial t} \wedge dt \wedge \psi \wedge d_/f \wedge f^*(\zeta) + d_/\xi \wedge \psi \wedge \frac{\partial f}{\partial t}.dt \wedge f^*(\zeta) $$
   puisque \ $f^*(ds) = df = d_/f +  \frac{\partial f}{\partial t}.dt $. On trouve alors, puisque \\ 
   $$\nabla(d_/\xi) = d_/f \wedge \frac{\partial \xi}{\partial t} -  \frac{\partial f}{\partial t}.d_/\xi  $$ 
   $$  < \frac{\partial F}{\partial t}.dt ,  ds \wedge \zeta > =  \int_{f^{-1}(\Delta)} \ \nabla(d_/\xi) \wedge \psi \wedge dt \wedge f^*(\zeta).$$
   On a donc bien, \`a nouveau gr\^ace au th\'eor\`eme de Fubini, la formule annonc\'ee, pour \ $d_/ \xi \in \mathbb{E}$ :
     $$ \frac{\partial F}{\partial t}(s,t) = \int_{\gamma_{s,t}} \ \frac{\nabla(d_/\xi)}{d_/f}. $$
     Remarquons que si \ $I : \mathbb{E} \to \mathcal{O}_{\mathbb{H}\times \mathbb{C}}$ \ est le morphisme de faisceau localement d\'efini par \ $[d_/\xi] \to \int_{\gamma_{s,t}} \ \frac{d_/\xi}{d_/f} $, on a \'etabli les \'egalit\'es suivantes :\\
     $$(\partial_s)_{\circ} I_{\circ} b = I \quad {\rm et} \quad (\partial_t)_{\circ} I = I_{\circ}(b^{-1}\nabla)  $$
     de morphismes de faisceaux respectivement de \ $\mathbb{E}$ \ \`a valeurs dans \ $\mathcal{O}_{\mathbb{H}\times \mathbb{C}}$ \ et de \ $\mathbb{P}$ \ \`a valeurs dans \ $\mathcal{O}_{\mathbb{H}\times \mathbb{C}}$. $\hfill \blacksquare $
     
     \bigskip

     Reprenons les notations introduites au d\'ebut de ce paragraphe. Le lemme suivant montre que les int\'egrales sur les cycles d\'eterminent une section du faisceau \ $\mathbb{E}$.
     
           \begin{lemma} 
     Soit \ $p$ \ un point de \ $S^*$ \ et soit \ $\theta$ \ un germe en \ $p$ \ de section du faisceau \ $\mathbb{E}$ \ tel que pour tout couple \ $(s,t)$ \ assez voisin de \ $(0, t(p))$ \ dans \ $\mathbb{C}^* \times \mathbb{C}$ \ on ait
     $$ \int_{\gamma^j_{s,t}} \ \frac{\theta}{d_/f} \  \equiv 0$$
     pour une base horizontale \ $\gamma^j_{s,t}$ \ d\'eduite d'une base de l'espace vectoriel   \\ 
     $$H_{n-1}(\{ f = \varepsilon\} \cap U, \mathbb{C})  \simeq H_{n-1}(\{f = \varepsilon\} \cap \{ t = t(p)\} \cap U, \mathbb{C}) .$$
      Alors on a \ $\theta = 0 $ \ dans \ $\mathbb{E}$.
     \end{lemma}
     
     \parag{Preuve} D'apr\`es notre hypoth\`ese, pour chaque \ $t_0$ \ fix\'e assez voisin de \ $t(p)$ \ la section du fibr\'e de Gauss-Manin de la fonction \ $f_{t_0}$ \ sur l'hyperplan \ $\{ t = t_0\}$, induite par \ $\theta_{\vert \{t = t_0 \}}$,  est de torsion. Comme la fonction \ $f_{t_0}$ \ a une singularit\'e isol\'ee en \ $p(t_0) : = S^* \cap t^{-1}(t_0) $, son fibr\'e de Gauss-Manin est sans torsion d'apr\`es [S. 70], et on a donc que la valeur de \ $\theta$ \ en \ $p(t_0)$ \ est nulle, c'est \`a dire que l'image de \ $\theta$ \ dans \ $\mathbb{E}/(t-t_0).\mathbb{E}$ \ est nulle. Comme ceci a lieu  pour tout \ $t_0$ \ assez voisin de \ $t(p)$ \ et comme le faisceau \ $\mathbb{E}$ \ est localement libre de rang fini sur \ $t^{-1}(\mathcal{O}_D[[b]])$ \ sur \ $S^*$, on en conclut que \ $\theta $ \ est nulle au voisinage de \ $p$. $\hfill \blacksquare$

\section{Le th\'eor\`eme de \ $b^{-1}.\nabla-$finitude sur \ $S^*$.}

\subsection{Preuve du th\'eor\`eme de finitude.}

Le th\'eor\`eme \ref{technique} sera une cons\'equence simple du r\'esultat suivant.

\begin{thm}\label{finitude} 
Dans la "situation consid\'er\'ee" pr\'ecis\'ee au d\'ebut du chapitre 2, \ $\mathbb{E}$ \ est \ $b^{-1}.\nabla-$fini sur \ $S^*$.
\end{thm}

Reprenons les notations du paragraphe 3.2.\\
Notons \ $q : \hat{\Omega}^{\bullet} \to  \hat{\Omega}_/^{\bullet}$ \ l'application quotient et remarquons que l'on a l'\'egalit\'e \ $ d_/\circ q = q \circ d $. Nous noterons encore par \ $d_/ : \hat{\Omega}^{\bullet} \to  \hat{\Omega}_/^{\bullet +1}$ \ cette application compos\'ee. Donc pour \ $\omega \in \hat{\Omega}^{n-1}$ \ nous aurons l'\'el\'ement \ $[ d_/\omega] \in \mathbb{E}$.\\
Si les formes \ $\omega, \omega' \in \hat{\Omega}^{n-1}$ \ v\'erifient
$$ f.d\omega = df \wedge \omega' $$
on aura dans \ $\mathbb{E}$ \ la relation \ $a.[d_/\omega] = b.[d_/\omega']$. Ceci se v\'erifie facilement en utilisant l'identification de \ $\Omega_/^{\bullet}$ \ avec les formes ne pr\'esentant pas l'\'el\'ement diff\'erentiel "$dt$" \ et la d\'ecomposition \ 
$$ \Omega^{\bullet} = \Omega_/^{\bullet} \oplus dt \wedge \Omega_/^{\bullet-1} .$$

  \begin{prop}
  Soit \ $p \in S^*$ \ et  consid\'erons une base de Jordan \ 
   $e_1, \cdots, e_k$ \ d'un bloc de Jordan de la monodromie  de \ $f$ \  pour la valeur propre \ $exp(2i\pi.u)$ \ o\`u \ $u \in [0,1[$, agissant sur \ $H^{n-1}(\{f  = \varepsilon \} \cap U, \mathbb{C})$. \\
   Alors il existe un entier  $m \geq 0$ \ et des $(n-1)-$formes holomorphes  \ $\omega_1, \cdots, \omega_k$ \ sur \ $U$ \ v\'erifiant les propri\'et\'es suivantes:
  \begin{enumerate}
  \item On a sur \ $U$ \ les relations \  $d\omega_j =  (m + u)\frac{df}{f}\wedge \omega_j + \frac{df}{f}\wedge \omega_{j-1} \quad \forall j \in [1,k]$ \ avec la convention \ $\omega_{0}\equiv 0 $.
  \item L'espace vectoriel engendr\'e par les classes de cohomologie induites sur \\
   $\{f  = \varepsilon \} \cap U$ \ par les \ $\omega_j$ \ est \'egal  au sous-espace vectoriel engendr\'e par \ $e_1, \cdots, e_k$ \ de \ $H^{n-1}(\{f  = \varepsilon \} \cap U, \mathbb{C})$.
      \end{enumerate}
   Dans ces conditions les sections de \ $\mathbb{E}$ \ induites par les formes \ $d_/\omega_j$ \ sont dans le sous-faisceau  \ $Ker\,\nabla_{\vert \Delta}$ \ et le sous-$\mathbb{C}[[b]]-$module qu'elles engendrent est libre de rang \ $k$ \ et stable par \ $a$. Sa fibre est le (a,b)-module \`a p\^ole simple de rang \ $k$ \ de \ $\mathbb{C}[[b]]-$base \ $\varepsilon_1, \cdots, \varepsilon_k$ \  avec \ $a.\varepsilon_j = (m+u).b.\varepsilon_j + b.\varepsilon_{j-1}$ \ avec la convention \ $\varepsilon_0 = 0$.\\
   Le sous-faisceau de \ $t^{-1}(\mathcal{O}_{\Delta}[[b]])-$module engendr\'e par \ $[d_/\omega_1], \cdots, [d_/\omega_k]$ \ est stable par \ $a$ \ et  libre de rang \ $k$ \ sur \ $t^{-1}(\mathcal{O}_{\Delta}[[b]])$. Il est stable par \ $b^{-1}.\nabla$ \ et donc contenu dans \ $\mathbb{G}$.\\
   Si on part d'une base de Jordan compl\`ete  de l'espace vectoriel   \ $H^{n-1}(\{f  = \varepsilon \} \cap U, \mathbb{C})$,  le sous-faisceau de \ $t^{-1}(\mathcal{O}_{\Delta}[[b]])-$modules ainsi obtenu est encore libre de rang fini.
   \end{prop}

  \parag{Remarque} Gr\^ace \`a notre trivialisation locale,  pour chaque \ $p' \in \Delta$, les restrictions de la base \ $e_1, \cdots, e_k$ \   \`a \ $\{f_{p'} = \varepsilon \}$ \ induiront une base de Jordan d'un bloc de Jordan de taille \ $k$ \  pour la monodromie de \ $f_{p'}$ \ agissant sur \ $H^{n-1}(\{f_{p'} = \varepsilon \}\cap U, \mathbb{C})$ \ (resp. une base de Jordan compl\`ete de \ $H^{n-1}(\{f_{p'} = \varepsilon \}\cap U, \mathbb{C})$ \ si on part d'une base de Jordan compl\`ete de  \ $H^{n-1}(\{f  = \varepsilon \} \cap U, \mathbb{C})$ ). $\hfill \square$
  
  \parag{Preuve} Commen{\c c}ons par remarquer que si \ $\omega \in \Omega^{n-1}$ \ on a
  $$ df \wedge d\omega = dt\wedge \nabla(d_/\omega) .$$
  Puisque l'on a \ $df \wedge d\omega_j = 0, j \in [1,k]$\ les classes \ $[d_/\omega_j] \in \mathbb{E}$ \ sont donc bien dans \ $Ker\,\nabla_{\vert \Delta}$.\\
    L'existence de l'entier \ $m$ \ et des $(n-1)-$formes holomorphes v\'erifiant les propri\'et\'es 1) et 2) r\'esulte de  [B.84]. Les sections sur l'ouvert \ $\Delta$ \ du faisceau \ $\hat{\mathcal{H}}^n$ \ induites par les formes \ $d_/\omega_j$ \ v\'erifient les relations :
  $$ a.[d_/\omega_j] = (m+u).b.[d_/\omega_j] + b.[d_/\omega_{j-1}] \quad \forall j \in [1,k]$$
  qui ne sont qu'une r\'e\'ecriture des relations 1).\\
   Montrons l'ind\'ependance sur  \ $\mathcal{O}_D[[b]]$. Si on a, sur un ouvert connexe \ $\Delta' \subset \Delta$, une relation 
  $$ \sum_{j = 1}^k \ s_j.[d_/\omega_j] \in d_/(\Gamma(\Delta', (\hat{K}er\, d_/f)^{n-1}) $$
  o\`u \ $s_j \in \Gamma(\Delta', \mathcal{O}_{D}[[b]])$, pour \ $j \in [1.k]$,  on obtient  pour chaque \ $t \in \Delta'$ \ une relation sur \ $\mathbb{C}[[b]]$ \ des classes correspondantes, ce qui  donne la nullit\'e des coefficients pour chaque \ $t \in \Delta'$ \ fix\'e. En effet, pour \ $t \in \Delta'$ \ fix\'e,  ces classes sont ind\'ependantes  sur  \ $\mathbb{C}\{f_t\}$ \ dans le syst\`eme de Gauss-Manin la fonction \ $f_t$ \ gr\^ace \`a la remarque qui suit l'\'enonc\'e de la proposition, et donc aussi dans le module de Brieskorn correspondant  puisqu'il est sans torsion d'apr\`es [S.70]. On en d\'eduit la m\^eme propri\'et\'e dans  son compl\'et\'e formel en \ $f_t$ \  par platitude de \ $\mathbb{C}[[f_t]]$ \ sur \ $\mathbb{C}\{f_t\}$, et celui-ci  co{\"i}ncide avec son compl\'et\'e formel en \ $b$, qui est le (a,b)-module associ\'e \`a \ $f_t$ \ \'egal \`a la fibre en \ $t$ \ de \ $\mathbb{E}$. L'ind\'ependance sur \ $\mathbb{C}[[b]]$ \ en r\'esulte.\\
  L'inclusion dans \ $\mathbb{G}$ \ est imm\'ediate puisque les \ $[d_/\omega_j]$ \ sont dans \ $Ker\,\nabla$ \ (voir prop. \ref{G} (2)).\\
  L'ind\'ependance pour une base  de Jordan  compl\`ete s'obtient de fa{\c c}on  analogue. $\hfill \blacksquare$

 \parag{D\'emonstration du th\'eor\`eme}  Compte tenu de la proposition pr\'ec\'edente, il nous suffit de montrer qu'il existe localement sur \ $S^*$ \ un entier \ $N$ \ tel que le  \ $t^{-1}(\mathcal{O}_{\Delta}[[b]])-$module \ $\tilde{\mathbb{G}} \subset \mathbb{G}$ \ construit \`a partir d'une base de Jordan de la monodromie de \ $f$ \ agissant sur l'espace vectoriel \  $H^{n-1}(\{f  = \varepsilon \} \cap U, \mathbb{C})$ \ v\'erifie \ $b^N.\mathbb{E} \subset \tilde{\mathbb{G}}$. Ceci r\'esulte du fait que si l'on consid\`ere la construction effectu\'ee, on peut choisir, d'apr\`es [B.86], l'entier \ $m$ \ au plus \'egal \`a \ $n$. On obtient alors, pour chaque singularit\'e isol\'ee \ $f_{p'}$, un r\'eseau \ $\tilde{\mathbb{G}}_{p'}$ \ du r\'eseau de Brieskorn \ $\mathbb{E}_{p'}$ \ qui v\'erifie \ $b^n.\mathbb{E}_{p'} \subset \tilde{\mathbb{G}}_{p'} $. Comme ceci a lieu pour chaque \ $p'$ \ assez voisin de \ $p$, on en conclut que \ $b^n.\mathbb{E} \subset \tilde{\mathbb{G}} \subset \mathbb{G}$. $\hfill \blacksquare$
  
  \subsection{Estimation de \ $\mathbb{G}$.}
  
  La preuve du th\'eor\`eme de \ $b^{-1}.\nabla-$finitude montre que \ $\mathbb{G}$ \ contient \  $b^n.\mathbb{E}$ sur \ $S^*$. Il est int\'eressant de pouvoir am\'eliorer cette "minoration" de \ $\mathbb{G}$ \ dans certains exemples. La proposition ci-dessous donne un tel crit\`ere.
  
  \bigskip
  
  Nous allons travailler pr\`es d'un point g\'en\'erique de \ $S^*$. Choisissons, pr\`es d'un tel point, un syst\`eme de coordonn\'ees locales \ $t, x_1, \dots, x_n$ \ tel que l'on ait 
  $$S^* = \{ x_1 = \cdots = x_n = 0\} .$$
   L'id\'eal \ $\frak{M}$ \ de \ $\mathcal{O}_X$ \ engendr\'e par \ $x_1, \dots, x_n$ \ est donc l'id\'eal r\'eduit de \ $S^*$.
  
   Commen{\c c}ons par remarquer que l'on a, localement sur \ $S^*$, \'equivalence entre l'\'egalit\'e  \ $\mathbb{G} = \mathbb{E}$ \ et l'appartenance de \ $\frac{\partial f}{\partial t}$ \ \`a \ $J_/(f)$, l'id\'eal jacobien relatif de \ $f$. En effet, d\`es que \ $\frac{\partial f}{\partial t} \not\in J_/(f)$, on a \ $dx \not\in \mathbb{P}$ \ et donc, \`a fortiori \ $dx \not\in \mathbb{G}$.
   
   \bigskip
   
  \begin{prop}\label{Estim.}
   Notons \ $\mathcal{M}^k$, pour \ $k \geq 1$, l'image de \ $\frak{M}^k.\hat{\Omega}^n_/$ \ dans \ $\mathbb{E}$.\\
    Supposons qu'il existe un entier \ $k $ \ tel que 
    \begin{equation*}
   \frak{M}^k.\frac{\partial f}{\partial t} \subset \frak{M}^{k+1}.J_/(f)  .\tag{*}
    \end{equation*}
  Alors le plus grand sous-$\mathcal{O}_D[[b]]-$module (coh\'erent)  \ $\mathbb{G}$ \ de \ $\mathbb{E}$, stable par \ $a$ \ et \ $b^{-1}\nabla$ \ contient \ $ \mathcal{M}^k$. 
  \end{prop}
    
    \bigskip
    
    \noindent \textit{Preuve.} Nous allons montrer que sous notre hypoth\`ese \ $G : = \mathcal{M}^k$ \  est stable par \ $a$ \ et \ $b$ \ et il v\'erifie \ $\nabla(G) \subset  b.G $.
   
  La stabilit\'e par \ $a$ \ est \'evidente. Pour montrer  la stabilit\'e par \ $b$, nous allons montrer que \ $b.G$ \ est l'image dans \ $\mathbb{E}$ \ de \ $\frak{M}^{k+1}.J_/(f).\Omega^n_/$, ce qui implique en particulier la stabilit\'e par \ $b$.\\
Si \ $\omega = d_/\xi$ \ avec \ $\omega \in \frak{M}^k.\Omega_/^n$, on peut choisir \ $\xi \in \frak{M}^{k+1}.\Omega_/^{n-1}$. On aura alors \ $ b[\omega] = [d_/f \wedge \xi] $ \ qui est bien dans \ $\frak{M}^{k + 1}.J_/(f).\Omega_/^n$.\\
    R\'eciproquement,  si \ $\eta \in \frak{M}^{k + 1}.J_/(f).\Omega_/^n$ \ on peut \'ecrire \ $\eta  = d_/f \wedge \xi $ \ avec \\
     $\xi \in \frak{M}^{k+1}.\Omega_/^{n-1}$. Alors \ $d_/\xi \in \frak{M}^k.\Omega_/^n$ \ et donne un \'element \ $[d_/\xi] \in G$ \ dont l'image par \ $b$ \ est \ $[\eta]$. Ceci prouve notre assertion.\\
    On a alors, pour \ $\omega = d_/\xi$,
    $$ \nabla(\omega) = d_/f \wedge \frac{\partial \xi}{\partial t} - \frac{\partial f}{\partial t}.\omega.$$
    Comme pour \ $[\omega] \in G$ \ on peut choisir \ $\xi$ \ et donc \ $ \frac{\partial \xi}{\partial t} \in \frak{M}^{k+1}.\Omega_/^{n-1}$, on obtient, gr\^ace \`a l'inclusion  (*), \
    $ \nabla(\omega) \in   \frak{M}^{k + 1}.J_/(f).\Omega_/^n \subset b.G. \hfill \blacksquare$

 \subsection{Exemples.}
 
  \parag{Exemple 1} On se propose d'\'etudier le cas o\`u \ $f(t,x) = P(x) + t.Q(x)$ \ avec   \ $P$ \ et \ $Q$ \ deux germes de fonctions holomorphes \`a l'origine de \ $\mathbb{C}^n$. \\
   
   \begin{lemma}
   Supposons \ $P$ \ \`a singularit\'e isol\'ee \`a l'origine et supposons que l'on ait \ $\frak{M}^{k+1}.J(Q) \subset \frak{M}^{k+1}. J(P)$ \ pour un entier \ $k \geq 0$, o\`u \ $\frak{M}$ \ d\'esigne l'id\'eal maximal de l'origine dans \ $\mathbb{C}^n$.\\
   Alors on a \ $\widehat{\frak{M}}^{k+1}.J_/(f) = \widehat{\frak{M}}^{k+1}.J(P)$ \ o\`u \ $\widehat{\frak{M}}$ \ d\'esigne l'id\'eal engendr\'e par \ $\frak{M}$ \ dans \ $\mathcal{O}_{\mathbb{C}^{n+1}}$.\\
   Si on a de plus \ $Q \in \frak{M}.J(Q)$, ce qui est v\'erifi\'e en particulier si \ $Q$ \ est quasi homog\`ene, on obtiendra l'inclusion \ $\frak{M}^{k}.\frac{\partial f}{\partial t} \subset \frak{M}^{k+1}.J_/(f)$ \ et la proposition pr\'ec\'edente donnera que \ $ \mathcal{M}^{k} \subset \mathbb{G}$ \ sur \ $S = \{ x = 0 \}$ \ au voisinage de l'origine.
   \end{lemma} 
   
   \parag{Preuve} L'hypoth\`ese permet donc d'\'ecrire chaque \ $x^{\alpha}.\frac{\partial Q}{\partial x_j}$ \ pour \ $\alpha \in \mathbb{N}^n$, v\'erifiant \ $ \vert \alpha \vert = k+1$, comme combinaison lin\'eaire \`a coefficients holomorphes des \ $x^{\beta}.\frac{\partial P}{\partial x_i}$. Si \ $\gamma$ \ d\'esigne le vecteur colonne des \ $x^{\alpha}.\frac{\partial P}{\partial x_j}$ \ et \ $\delta$ \ le vecteur colonne des \ $x^{\alpha}.\frac{\partial Q}{\partial x_j}$ \ on aura donc \ $\delta = \mathcal{R}.\gamma$ \ o\`u \ $\mathcal{R}$ \ est une matrice \`a coefficients holomorphes dans \ $\mathbb{C}^n$.\\
 Comme les \ $x^{\alpha}.\big[\frac{\partial P}{\partial x_j} + t.\frac{\partial Q}{\partial x_j}\big] $ \ forment un syst\`eme g\'en\'erateur sur \ $\mathcal{O}_{\mathbb{C}^{n+1}}$ \ de \ $ \widehat{\frak{M}}^{k+1}.J_/(f) $, on aura une relation matricielle
   $$ \Gamma = (Id +t.\mathcal{R}).\gamma $$
   o\`u \ $\Gamma$ \ d\'esigne le vecteur colonne des \ $x^{\alpha}.\big[\frac{\partial P}{\partial x_j} + t.\frac{\partial Q}{\partial x_j}\big] $. Pour \ $\vert t \vert \ll 1 $ \ la matrice \ $Id + t. \mathcal{R}$ \ sera inversible et on a ainsi \'etabli l'\'egalit\'e \ $\widehat{\frak{M}}^{k+1}.J_/(f) = \widehat{\frak{M}}^{k+1}.J(P)$ \ au voisinage de l'origine dans \ $\mathbb{C}^{n+1}$. Cette \'egalit\'e montre que le lieu singulier de \ $f$ \ est contenu dans \ $S = \{ x = 0 \}$ \ au voisinage de l'origine, et comme on a \ $Q(0)=0$, \'egalit\'e impos\'ee par l'appartenance de \ $Q$ \ \`a \ $\frak{M}.J(Q)$, on aura \'egalit\'e. On conclut imm\'ediatement car \ $Q \in J_/(f)$ \ implique \ $Q \in J(P)$. $\hfill \blacksquare$
   
    \bigskip
    
       Illustrons ceci par un exemple simple ''explicite''.

  \parag{Exemple 2}   Il s'agit de calculer l'exemple suivant dans lequel on a  \ $n = 2$ : \\ 
  $P(x,y) = x^4 + y^4, Q(x,y) = x^2.y^2 $ \ et donc \ $f(x,y,t) = P(x,y) + t.Q(x,y)$ \ qui est une d\'eformation \`a \ $\mu-$constant  pour \ $t \not= \pm 2$, d'hypersurfaces \`a singularit\'es isol\'ees de \ $\mathbb{C}^2$ \ dont le faisceau des modules de Brieskorn n'est pas   localement constant. En effet le birapport des quatre droites de \ $\mathbb{C}^2$ \ que l'on obtient pour chaque valeur de \ $t$ \ et qui vaut\footnote{Il s'agit de calculer le birapport des quatre racines de l'\'equation \ $z^4 +t.z^2 + 1 = 0 $. Ce nombre est d\'efini modulo le groupe du birapport. Par exemple \ $1 - (-\frac{t - 2}{4}) = \frac{t+2}{4}$ \ repr\'esente la m\^eme classe de birapport.} \  \ $-\frac{t - 2}{4}$ \ est une fonction localement injective, ce qui montrent que ces germes de fonctions holomorphes ne sont jamais localement deux \`a deux analytiquement \'equivalentes.\\
Cependant l'hypoth\`ese \ $\frak{M}^2.J(Q) \subset \frak{M}^2.J(P) $ \ du lemme pr\'ec\'edent est satisfaite, o\`u ici, on a simplement \ $\frak{M} : = (x,y)$. \\
On v\'erifie aussi imm\'ediatement  que \ $Q \not\in J(P)$. On aura donc \ $\mathbb{G} = \mathcal{M}$.

\bigskip

Explicitons le calcul de \ $\nabla$.

On a 
\begin{align*}
\frac{\partial f}{\partial x} = 4x^3 + 2t.xy^2 \\
\frac{\partial f}{\partial y} = 4y^3 + 2t.x^2y.
\end{align*}
et donc
\begin{align*}
8.x^3y & = 2y \frac{\partial f}{\partial x} - tx.(\frac{\partial f}{\partial y} - 2t.x^2y) \\
\quad   & = 2t^2.x^3y +  2y \frac{\partial f}{\partial x} - tx.\frac{\partial f}{\partial y}.
\end{align*}
On obtient ainsi
\begin{align}
2(4 - t^2).x^3y & =  2y \frac{\partial f}{\partial x} - tx.\frac{\partial f}{\partial y} \in J_/(f) \quad {\rm et}\\
2(4 - t^2).xy^3 & =  2x \frac{\partial f}{\partial y} - ty.\frac{\partial f}{\partial x} \in J_/(f) \quad {\rm par \ symetrie}.
\end{align}
On en d\'eduit que
\begin{align}
4.x^5 & = x^2.\frac{\partial f}{\partial x} - 2t.x^3y^2 \in \in J_/(f) \quad {\rm et}\\
4.y^5 &  = y^2.\frac{\partial f}{\partial y} - 2t.x^2y^3 \in J_/(f) \quad {\rm par \ symetrie}.
\end{align}
Gr\^ace aux relations
\begin{align}
4.x^4 = - 2t.x^2y^2 + x.\frac{\partial f}{\partial x} \\
4.y^4 = -2t.x^2y^2 + y.\frac{\partial f}{\partial y}
\end{align}
une \ $\mathcal{O}_{D}[[b]]-$base de \ $\mathbb{E}$ \ est donc donn\'ee par
$$ (1, x, y, x^2, y^2, xy, x^2y, xy^2, x^2y^2).dx\wedge dy.$$
Soit \ $\omega : = x.dy - y.dx $. Pour un mon\^ome \ $m$ \ homog\`ene de degr\'e \ $\delta(m)$ \ on a
\begin{equation}
d_/(m.\omega) = \frac{ \delta(m) + 2}{4}.\frac{d_/f}{f}\wedge m.\omega 
\end{equation}
ce qui donne 
\begin{equation}
a(m) =  \frac{ \delta(m) + 2}{4}.b(m).
\end{equation}
On a \'egalement
\begin{align}
\nabla(d(m.\omega)) = d_/f\wedge \frac{\partial(m.\omega)}{\partial t} - x^2y^2.d(m.\omega)\\
{\rm et \ donc} \quad  \nabla(m) = -x^2y^2.m.
\end{align}
Comme \ $x^2y^2.\frak{M} \subset J_/(f)$, on voit que \ $b^{-1}\nabla$ \ op\`ere sur \ $\mathbb{G} = \mathcal{M}$.\\
Par exemple, comme
\begin{align} 
2(4 - t^2).x^3y^2 & = 2y^2 \frac{\partial f}{\partial x} - txy.\frac{\partial f}{\partial y}\\
\quad                     & = d_/f \wedge (2y^2dy + t.xydx)
\end{align}
on aura 
$$b^{-1}\nabla(x) = \frac{t.x}{2(4 -t^2)}.$$
Donc \ $(4 -t^2)^{\frac{1}{4}}.x$ \ sera (localement)  dans \ $Ker\, \nabla$ \ pour \ $\vert t \vert < 2$.\\

\smallskip

Comme \ $\nabla(1) = - x^2y^2 \notin b.\mathbb{E}$, pour tester directement  la \ $b^{-1}\nabla-$finitude  de \ $\mathbb{E}$, posons \ $E_1 : = \mathbb{E}\oplus \mathcal{O}_{D}.\varepsilon $, o\`u l'on d\'efinit \ $\varepsilon : = b^{-1}(x^2y^2)$. \\
Alors on obtient, puisque \ $a(x^2y^2) = \frac{3}{2}b(x^2y^2)$ \ les relations suivantes :
$$  a\varepsilon = \frac{1}{2}.b\varepsilon \quad {\rm et} \quad \nabla{\varepsilon} = b^{-1}\nabla(x^2y^2) = b^{-1}(-x^4y^4).$$
Explicitons \ $ b^{-1}(-x^4y^4)$. On a d'apr\`es (1)
\begin{equation}
2(4 - t^2).x^4y^4 =  d_/f \wedge (2xy^4dy + t.x^2y^3dx)
\end{equation}
ce qui donne \ $ 2(4 - t^2).x^4y^4  = b(2y^4 - 3t.x^2y^2)$ \ et donc d'apr\`es (6)
$$ 2(4 - t^2).x^4y^4 = b( -4t.x^2y^2 + d_/f \wedge (-\frac{1}{2}y.dx)) = b( -4t.x^2y^2 + \frac{1}{2}.b(1)) .$$
On a donc 
 $$ b^{-1}\nabla(\varepsilon) =  \frac{2t}{4 - t^2}.\varepsilon  - \frac{1}{4(4-t^2)}.1 .$$
Ceci permet de conclure que \ $E_1$ \ est stable par \ $b^{-1}\nabla$. \\
On obtient ainsi directement la \ $b^{-1}\nabla-$r\'egularit\'e pour \ $\mathbb{E}$ \ dans cet exemple.

\smallskip

Le \ $\mathcal{O}_D[[b]]-$module \ $\mathcal{O}_D[[b]].1 \oplus  \mathcal{O}_D[[b]].\varepsilon $ \ est stable par \ $a$ \ et \ $b^{-1}\nabla$, via les formules:
\begin{align*}
& a.\varepsilon = \frac{1}{2}.b\varepsilon , \quad a.1 =  \frac{1}{2}.b.1 \\
& b^{-1}\nabla(1) = -\varepsilon , \quad b^{-1}\nabla(\varepsilon) = \frac{2t}{4 - t^2}.\varepsilon  - \frac{1}{4(4-t^2)}.1 .
\end{align*}
    
   \parag{Exemple 3}    
    Soient \ $p,q,r$ \ trois entiers \ $\geq 3 $ \ tels que \ $\frac{1}{p} + \frac{1}{q} + \frac{1}{r} < 1 $ \ et posons
    $$ f(x,y,z,t) = x^p + y^q + z^r + txyz .$$
    Remarquons que l'on a \ $f \in \frak{M}^3$ \ et donc \ $J_/(f) \subset  \frak{M}^2$. \\
    Montrons que, sur l'ouvert \ $S^* : = \{ t \not= 0 \} $ \ de \ $\mathbb{C}$, les hypoth\`eses de la proposition \ref{Estim.} sont v\'erifi\'ees avec \ $k = 1$. D'abord on a
    $$ J_/(f) = (p.x^{p-1} + t.yz, q.y^{q-1} + t.xz, r.z^{r-1} + t.xy) .$$
    Les relations 
    $$\frac{\partial f}{\partial x} = p.x^{p-1} + t.yz \quad \frac{\partial f}{\partial y} = q.y^{q-1} + t.xz \quad \frac{\partial f}{\partial z} = r.z^{r-1} + t.xy $$
    donnent
    \begin{align*}
    pq.x^{p-1}.y^{q-1} = t^2.xyz^2 + \frak{M}^2.J_/(f) \quad {\rm ainsi \ que} \\
    t.x^{p-1}.y^{q-1} = - r.z^{r-1}.x^{p-2}y^{q-2} + \frak{M}^2.J_/(f) 
    \end{align*}
    puisque \ $x^{p-2}y^{q-2} \in \frak{M}^2$. On a donc
    $$ t^3.xyz^2 + pqr.z^{r-1}.x^{p-2}y^{q-2} + \frak{M}^2.J_/(f) $$
    ou encore, puisque \ $p, q, r, $ \ sont au moins \'egaux \`a $3$ \ et \ $t\not=  0$,
    \begin{equation*}
     t^3.xyz^2(1 + \frac{pqr}{t^3}.z^{r-3}x^{p-3}y^{q-3}) \in \frak{M}^2.J_/(f) . \tag{@}
     \end{equation*}
     On aura aussi
     \begin{align*}
     & t.x^py  & = -r. x^{p-1}z^{r-1} + x^{p-1}.\frac{\partial f}{\partial z} \\
     &t^2.x^py  & = -r.t.x^{p-1}z^{r-1} +  \frak{M}^2.J_/(f) = rq.y^{q-1}.x^{p-2}.z^{r-2} +  \frak{M}^2.J_/(f) \\
     & \qquad & = rq.y^2.x.z.(y^{q-3}.x^{p-3}.z^{r-3}) +  \frak{M}^2.J_/(f) \in  \frak{M}^2.J_/(f) 
     \end{align*}
     d'apr\`es (*).
    On a donc, en utilisant encore le fait que \ $x, y, z$ \ jouent le m\^eme r\^ole, que
    $$\frak{M}.\frac{\partial f}{\partial t} \subset\frak{M}^2.J_/(f) \quad {\rm ainsi \ que} \quad \frak{M}.f \subset \frak{M}^2.J_/(f) .$$
    Il nous reste seulement \`a voir que \ $\frac{\partial f}{\partial t}.h \in J_/(f) $ \ implique \ $h \in \frak{M}$; ceci r\'esulte du fait que \ $\frac{\partial f}{\partial t} \not\in J_/(f) $ .\\
    
    \begin{lemma}
    Notons par \ $G : = \frak{M}.\Omega^n_/\big/ d_/f \wedge d_/\Omega^{n-2}_/ $.\\
    Comme le fibr\'e vectoriel sur \ $S^*$ \ d\'efini par \ $\mathbb{E}\big/ b.\mathbb{E} $ \ est trivial, on constate facilement que le \ $\mathcal{O}_S[[b]]-$module 
    $$\Gamma : = G + \mathcal{O}_S[[b]]. b^{-1}\alpha$$
     de \ $\mathbb{E}[b^{-1}]$ \ o\`u \ $\alpha : = xyz$, est stable par \ $a, b, b^{-1}\nabla$ \ et \ $b^{-1}a$ \ sur \ $S^*$.
    \end{lemma}
    
    \bigskip
    
    \noindent \textit{Preuve.} Comme on a  
    $$ p.x^p = q.y^q = r.z^r = - t.\alpha \quad modulo \quad\frak{M}.J_/(f) ,$$
    on aura
    \begin{equation}\label{@@@} 
     f - (1 - \rho).t\alpha \in \frak{M}.J_/(f)\quad{\rm avec} \quad \rho = \frac{1}{p} + \frac{1}{q} + \frac{1}{r}.
     \end{equation}
     
 On en d\'eduit en particulier que \ $a.\alpha \in \frak{M}^2.J_/(f)$ \ et donc que \ $ab^{-1}\alpha $ \ est bien \\ dans \ $\Gamma$.
 La stabilit\'e par \ $b^{-1}a$ \ en d\'ecoule alors puisque la relation \ $\frak{M}.f \subset \frak{M}^2.J_/(f)$ \   donne 
 $$a.G \subset \frak{M}^2.J_/(f) \subset b.G .$$
 
 \bigskip
 
\parag{Remarques}
\begin{enumerate}
\item  Le cas \ $ p = q = r  = 3 $ \ rel\`eve du premier exemple trait\'e plus haut. 
\item  Les cas  \ $\frac{1}{p} + \frac{1}{q} + \frac{1}{r} < 1 $ \ r\'esulte de l'\'etude de [B.II] puisque l'on a
 $$ W : = \frac{1}{p}.x.\frac{\partial}{\partial x} + \frac{1}{q}.y.\frac{\partial}{\partial y} + \frac{1}{r}.z.\frac{\partial}{\partial z}  + (1-\rho).t.\frac{\partial}{\partial t} $$
 qui v\'erifie \ $W.f = f$ \ et ne s'annule pas sur \ $\{ t \not= 0 \}$.
\item La relation   \ $\frak{M}.f \subset \frak{M}^2.J_/(f)$ \   implique que le plus grand sous-(a,b)-module \`a p\^ole simple de \ $f_t$ \ pour chaque \ $t \not= 0 $ \ est la fibre en \ $t$ \ de \ $\mathcal{M}$. \\
On a donc \ $\mathbb{G} = \mathcal{M} = \frak{M}.\mathbb{E}$ \ dans ce cas. Donc \ $\mathcal{M}_{t = 1}$ \ est le plus grand sous-(a,b)-module \`a pole simple de \ $\mathbb{E}_{ t=1}$, d'apr\`es la description de \ $\mathcal{H}^n$ \ dans le cas o\`u il existe un champ de vecteur holomorphe ne s'annulant pas et v\'erifiant \ $W.f = f$ \ pr\`es des points  de \ $S^*$ \ (voir [B.II] th\'eor\`eme 4.3.1.) 
\end{enumerate}

\newpage

    \section{R\'ef\'erences.}
     
     \begin{enumerate}
     
     \item{[B.84]} Barlet, D. \textit{Contribution effective de la monodromie aux d\'eveloppements asymptotiques,} Ann. Sc. Ec. Norm. Sup. 17 (1984), p.293-315.
     
      \item{[B.86]}  Barlet, D. \textit{Monodromie et p\^oles de \ $\int_X \vert f \vert^{2\lambda}\square$.} Bull. Soc. Math. France, t.114 (1986), p. 247-269.

  \item{[B.95]} Barlet, Daniel \textit{Theorie des (a,b)-modules II. Extensions} in Complex Analysis and Geometry, Pitman Research Notes in Mathematics Series  366 Longman (1997), p. 19-59.

 \item{[B.05]} Barlet, D. \textit{ Modules de Brieskorn et formes hermitiennes pour une singularit\'e isol\'ee d'hypersurface,}  Revue  de l'Institut E. Cartan (Nancy) vol.18 (2005),  p.19-46.

 \item{[B.II)]} Barlet, D. \textit{Sur certaines singularit\'es non isol\'ees d'hypersurfaces II}, \`a para\^itre au Journal of Algebraic Geometry.
  
  \item{[M.74]} B.Malgrange : {\it Int\'egrale asymptotique et monodromie.} Ann. Scient. Ec. Norm. Sup. , t.7 (1974), p.405-430
  
  \item{[S.70]} M. S\'ebastiani : {\it Preuve d'une conjecture de Brieskorn.} Manuscripta Math. 2 (1970), p.301-308.

\end{enumerate}

\bigskip

Barlet Daniel, Institut Elie Cartan UMR 7502  \\
Nancy-Universit\'e, CNRS, INRIA  et  Institut Universitaire de France, \\
BP 239 - F - 54506 Vandoeuvre-l\`es-Nancy Cedex.France.\\
e-mail : barlet@iecn.u-nancy.fr.

\end{document}